*Draft of June 2023—please do not cite without permission.*

# Montague's Paradox without Necessitation[1]

T. Parent (Nazarbayev University)
nontology@gmail.com

## *1. Introduction*

Standardly, 'necessarily' is treated in modal logic as an operator on propositions (much like '~'). Some have thought it should be seen instead as a predicate 'N(*x*)' defined on sentences. If that is right, then in lieu of a formula "□*p*," one would have something along the lines of "N('*p*')". But the operator view has become standard, mainly because of Montague's (1968) argument that liar-like paradox results under the predicate view, for systems T and stronger (see Halbach 2009). However, Dean (2014) argues that Montague's paradox relies on contentious use of the necessitation rule; his remedy, accordingly, is to restrict the necessitation rule. Yet the following demonstrates that Montague-like paradox also results without the use of the necessitation rule. This suggests not only that Dean's proposal is lacking, but more broadly, that any restriction on the necessitation rule will be insufficient to prevent this sort of paradox.

## *2. Montague's version of the paradox*

Montague originally presented his paradox not within modal logic *per se*, but rather in a version of Q (Robinson's 1950 arithmetic), where the system is assumed to contain a provability-predicate '*P*' governed by the following axiom-scheme and rule of inference. If $\ulcorner\phi\urcorner$ is the Gödel number of a formula $\phi$:[2]

[1] Thanks to Michael Detlefsen, Volker Halbach, and especially Alexander Pruss and for helpful feedback on earlier versions of this paper. Pruss (independently) offers a similar argument on his blog; see http://alexanderpruss.blogspot.com/2009/08/modal-liar.html.

($T_P$) $P(\ulcorner\phi\urcorner) \supset \phi$

($N_P$) If $\vdash\phi$, then $\vdash P(\ulcorner\phi\urcorner)$

Montague's paradox is then shown as follows:

| | | |
|---|---|---|
| (i) | $\vdash \Delta \equiv {\sim}P(\ulcorner\Delta\urcorner)$ | [From Gödel's Diagonal Lemma] |
| (ii) | $\vdash P(\ulcorner\Delta\urcorner) \supset \Delta$ | [From ($T_P$)] |
| (iii) | $\vdash {\sim}P(\ulcorner\Delta\urcorner)$ | [From (ii) and (i) by truth-functional logic] |
| (iv) | $\vdash \Delta$ | [From (i) and (iii)] |
| (v) | $\vdash P(\ulcorner\Delta\urcorner)$ | [From (iv) by ($N_P$)] |
| (vi) | $\vdash \bot$ | [Contradiction at (iii) and (v)] |

This can be seen as another illustration of how a scheme like ($T_P$) creates inconsistency in Q; Löb's theorem is another route for making the point.[3] Thus, such a scheme is typically rejected in formal arithmetic.[4]

But if $P(\ulcorner\phi\urcorner)$ is interpreted to mean that the *truth* of $\phi$ is *logically necessary*, ($T_P$) would be widely accepted. Concurrently, when the propositional operators '$\Box$' and '$\Diamond$' are used to express alethic modal notions, the following is commonly accepted as an axiom:

---

[2] Here and elsewhere, I gloss the use/mention distinction. Quine's (1951) quasi-quotes could be deployed to distinguish the formula $\ulcorner\phi\urcorner$ from the proposition $\phi$ expressed by $\ulcorner\phi\urcorner$. (Some other notation for Gödel numbers would then be required.) But to avoid clutter, I instead rely on context to disambiguate.
[3] See Boolos et al. (2007), ch. 18.
[4] Post (1970) also offers a modal liar paradox akin to Montague's, which also appears to assume a kind of necessitation. ("If a statement entails its negation, then its negation is necessary;" p. 405.) Post regards his paradox as different from Montagues', however, in that Post utilizes a *sentential prefix* 'It is possibly false that…' composing with a sentence to form a longer sentence. In contrast, Montague's uses a predicate '*x* is possibly false', where '*x*' is replaced by a singular term for a sentence (or rather, its Gödel number.) Regardless, the main issue in this paper is whether necessitation is needed for modal paradox; Post's and Montague's paradoxes do not differ in this regard.

(T) $\Box\phi \supset \phi$

Notwithstanding, when modality is expressed with operators, we cannot use the Diagonal Lemma to show that there is sentence which is modally equivalent to its own non-necessity. For the Lemma applies only in relation to predicates. Accordingly, the counterpart of (i) is not a theorem in such logics, and so (T) does not give rise to Montague's Paradox in alethic modal logics.

Be that as it may, Dean (2014) suggests that the proper response to Montague's Paradox is rather to reject or at least restrict the analogue of ($N_P$), a.k.a. the necessitation rule:

(N) If $\vdash\phi$, then $\vdash\Box\phi$

It is unclear whether Dean would also have us return to symbolizing modal notions with predicates rather than operators. For familiarity's sake, however, I shall continue to use the operators. My main purpose is to show that, contra Dean, even the complete absence of the necessitation rule would not alone suffice to block Montague's paradox. Or at least, that is so on the hypothesis that some sentence is equivalent to the negation of its own necessitation.

I shall in fact give two arguments to this effect. The first derives the paradox without necessitation in a variant of modal System T. However, since the derivation suggests that the necessitation rule is redundant in the system, the paradox in this case might not be unexpected. Thus, I also demonstrate that the paradox occurs in a weaker system which is not even complete with respect to truth-functional validities, much less modal validities from System T. This will make clear that the paradox depends in no way on necessitation.

*3. The paradox and system T*

We shall be using a natural deduction formulation of modal system T provided by Garson (2013) (albeit with a few superficial differences). Well-formed formulae (wffs) are defined in system T inductively as follows:

(1) $\bot$ and each propositional constant A, B, C,…$A_1$, $B_1$, $C_1$,…$A_2$, $B_2$, $C_2$,…are wff.

(2) If $\phi$ is a wff, then so are $\sim\phi$, $\Diamond\phi$, and $\Box\phi$.

(3) If $\phi$ and $\psi$ are wffs, then so are $(\phi \supset \psi)$ and $(\phi \equiv \psi)$.

(4) Nothing else is a wff.

(Per usual, the outermost parentheses are dropped when there is no danger of confusion.) $\bot$ is always assigned the truth-value False; the semantics for other wff is determined as follows (cf. Kripke 1959; 1963). Let a world *w* include a complete set of pairings $\langle\alpha, V\rangle$, where $\alpha$ is a propositional constant and V is the truth-value True or the truth-value False. A *complete* set of such pairings is such that, for every propositional constant $\alpha$, $\langle\alpha, \text{True}\rangle \in w$ or $\langle\alpha, \text{False}\rangle \in w$, but not both. Let *w*A*v* express that *w* bears the accessibility-relation to a world *v*. Then, the truth values for non-atomic wff in *w* are determined as follows:

(5) $\sim\phi$ is True in *w* iff $\phi$ is False in *w*.

(6) $(\phi \supset \psi)$ is True in *w* iff: $\phi$ is False or $\psi$ is True.

(7) $(\phi \equiv \psi)$ is True in *w* iff: $\phi$ and $\phi$ are both True or both False.

(8) $\Diamond\phi$ is True in *w* iff: There is a world *v* such that *wAv* and $\phi$ is True in *v*.

(9) $\Box\phi$ is True in *w* iff: For every world *v* such that *wAv*, $\phi$ is True in *v*.

The inference rules for system T include rules for introducing and eliminating $\bot$ and the operators. For $\bot$ and the truth-functional operators, the rules are as follows.

$\bot$ Intro

$\phi$

$\sim\phi$

$\therefore \bot$

$\bot$ Elim

$\bot$

$\therefore \phi$

$\sim$ Intro

| $\phi$
| $\vdots$
| $\bot$

$\therefore \sim\phi$

$\sim$ Elim

| $\sim\phi$
| $\vdots$
| $\bot$

$\therefore \phi$

$\supset$ Intro

| $\phi$
| $\vdots$
| $\psi$

$\therefore \phi \supset \psi$

$\supset$ Elim

$\phi \supset \psi$

$\phi$

$\therefore \psi$

$\equiv$ Intro

| $\phi$
| $\vdots$
| $\psi$

| $\psi$
| $\vdots$
| $\phi$

$\therefore \phi \equiv \psi$

$\equiv$ Elim

| | | |
|---|---|---|
| $\phi \equiv \psi$ | *-or-* | $\phi \equiv \psi$ |
| $\phi$ | | $\psi$ |
| $\therefore \psi$ | | $\therefore \phi$ |

For modal operators, the rules of inference are:

□ Elim

□ϕ

∴ ϕ

◇ Intro

ϕ

∴ ◇ϕ

□ Intro

| □
|———
| ⋮
| ϕ

∴ □ϕ

◇ Elim

◇ϕ

| □
|———
| | ϕ
| |———
| | ⋮
| | ψ

∴ ◇ψ

Modal Interchange Rules #1–#4

~□ϕ

∴ ◇~ϕ

◇~ϕ

∴ ~□ϕ

□~ϕ

∴ ~◇ϕ

~◇ϕ

∴ □~ϕ

A few remarks. First, regarding □ Intro, the subproof starting with '□' indicates what holds in an arbitrary accessible world. The idea behind the rule, then, is that if we can show that ϕ holds in an arbitrary world, then ϕ must hold of necessity. Thus, □ Intro is our system's version of the necessitation rule, and I shall often refer to it as such.

As concerns ◇ Elim, the subproof starting with '□' similarly indicates that we are considering what holds in an arbitrary accessible world. The ensuing sub-subroof then indicates a provisional assumption that ϕ holds in that world (where ◇ϕ is antecedently established). The rational here is that if it is known that ϕ is possible, then if ψ is shown in some arbitrary world where ϕ is assumed to hold, we thereby know that ψ is also possible.

Finally, note that □ Elim allows one to prove what one would otherwise be provable via the T-axiom, whereas ◇ Intro allows one to derive what would otherwise be derivable by the law of *ab esse ad posse*. There is therefore redundancy in the rules, due to the fact that *ab esse ad posse* is equivalent to the T-axiom, but no harm done. Similarly, I also include the following (redundant) rule for later ease of exposition:

<u>Impossibility of ⊥</u>

∴ ~◇⊥

The rule allows one to enter ~◇⊥ at any line; this just reflects that in no world is the False assigned the truth-value "True."

The aim, again, is to show that Montague's paradox can be generated without use of the necessitation rule (□ Intro). As mentioned earlier, however, the Diagonal Lemma cannot be used in the present context to show that some sentence of the language is modally equivalent to its own negated necessitation, which is a crucial ingredient to the paradox. So instead, I shall follow Dean (2014) and hypothesize that there is such a sentence μ, from which we will draw out an absurdity.

| | | |
|---|---|---|
| 1. | $\Box(\mu \equiv \sim\Box\mu)$ | Hypothesis |
| 2. | $\mu \equiv \sim\Box\mu$ | 1, $\Box$ Elim |
| 3. | $\Box\mu$ | [Assume for ~Intro] |
| 4. | $\mu$ | 3, $\Box$ Elim |
| 5. | $\sim\Box\mu$ | 2, 4, $\equiv$ Elim |
| 6. | $\bot$ | 3, 5, $\bot$ Intro |
| 7. | $\sim\Box\mu$ | 3-6, ~Intro |
| 8. | $\Diamond\sim\mu$ | 8, Interchange Rule #1 |
| 9. | $\Box$ | [Assume for $\Diamond$ Elim] |
| 10. | $\sim\mu$ | [Assume for $\Diamond$ Elim] |
| 11. | $\mu \equiv \sim\Box\mu$ | 1, $\Box$ Elim |
| 12. | $\Box\mu$ | [Assume for ~ Intro] |
| 13. | $\mu$ | 12, $\Box$ Elim |
| 14. | $\sim\Box\mu$ | 11, 13, $\equiv$ Elim |
| 15. | $\bot$ | 12, 14, $\bot$ Intro |
| 16. | $\sim\Box\mu$ | 12-15, ~ Intro |
| 17. | $\mu$ | 11, 16, $\equiv$ Elim |
| 18. | $\bot$ | 10, 17, $\bot$ Intro |
| 20. | $\Diamond\bot$ | 8, 9-18, $\Diamond$ Elim |
| 21. | $\sim\Diamond\bot$ | Impossibility of $\bot$ |
| 22. | $\bot$ | 20, 21, $\bot$ Intro |

Thus, it is refutable that some sentence $\mu$ is modally equivalent to its own negated necessitation. And at no point is the necessitation rule deployed.

Yet upon further inspection, the existence of this proof should not be too surprising. Rrather than using the necessitation rule to prove $\Box\mu$ directly under the hypothesis, we prove it indirectly by *reductio*. And the *reductio* proceeds by simply repeating (in relation to some accessible world) the earlier derivation of $\sim\Box\mu$. (Compare steps 12-16 with steps 3-7.) We are then able to derive $\mu$. However, since $\sim\mu$ is also assumed to hold in the relevant world, a contradiction follows. Therewith, we complete a necessitation-free derivation of $\Box\mu$.

Granted, if we did not have the rule "Impossibility of ⊥," we would instead need to derive $\Box\sim\bot$ using □ Intro in order to contradict $\Diamond\bot$ at line 20. But the Impossibility of ⊥ is much less general (hence, much weaker) than □ Intro. Even granting that, however, our indirect, necessitation-free proof-strategy seems generalizable.[5] Whenever $\phi$ is a theorem, we can assume that $\sim\phi$ holds in some accessible world, and then repeat the proof of $\phi$ in that world. We can thereby prove $\Box\phi$ by *reductio* without using the necessitation rule. Thus, if □ Intro were removed from the system, we apparently would not get a weaker system. And for our purposes, that is concerning. It now seems that our derivation of Montague's paradox is not independent of the necessitation rule, even if we do not invoke such a rule explicitly.

*4. The paradox in system T-minus*

For this reason, it is desirable to show that the paradox can also be derived without the necessitation rule in a patently weaker system. And in fact, the derivation above suggests how this can be done. Consider the system T-minus ("T–") which is just like system T, except that its inference rules are much more scant. The modal rules of inference are exhausted by □ Elim, ◇ Elim, and the *first* modal interchange rule. The only other inference rules in T– are ⊥ Intro, ~ Intro, and the Impossibility of ⊥.

[5] I thank Alexander Pruss for drawing my attention to this point.

The absence of ~ Elim (and the double negation rule) means that T– is incomplete with respect to truth-functional validities. But T– is also missing some modal theorems from T, given its limited stock of modal rules. One way to see the modal incompleteness is in trying to prove □~(A ≡ ~A) in T–. Since □ Intro is not part of T–, the only way to prove this would be by *reductio*, i.e., by assuming ~□~(A ≡ ~A) and deriving a contradiction. Yet since the system has does not have the rule ~ Elim (nor double negation), the *reductio* would at most allow us to conclude ~~□~(A ≡ ~A) via ~ Intro. So □~(A ≡ ~A) is an example of something which cannot be proved in T–, even though it can be proved in T.

Nevertheless, we can still prove Montague's paradox in T–. Indeed, the derivation of the paradox is the same as in T. Thereby, it is revealed that ⊥ is a consequence of line 1, and the hypothesis is refuted.

The proof in T– has the same strategy of "proof repetition" inside an accessible world. But unlike in T, that strategy does not generalize. Again, there are theorems like ~(A ≡ ~A) which cannot be necessitated by means of a *reductio* in T–. So the point stands that T– in no sense presupposes the necessitation rule.

*5. A point about novelty*

One colleague, who prefers anonymity, suggests that the above proof illustrates nothing new, for it is known that contradiction results whenever a system contains an operator T that conforms to the following conditions, where ϕ is any wff:

(C1) From Tϕ, ϕ is derivable.

(C2) From ~Tϕ, ~ϕ is derivable.

(See Friedman & Sheard 1987.) Note that when T is the box operator, $\Box$ Elim is sufficient for (C1) to hold of the system T–. And the rule of $\Diamond$ Elim in combination with ~ Intro is apparently sufficient for (C2) to hold. After all, given $\sim\Box\mu$ at line 7 (above), we show a contradiction using $\Diamond$ Elim. And given a contradiction, it is derivable using ~ Intro that $\sim\mu$. Such a derivation from $\sim\Box\mu$ to $\sim\mu$ is illustrative of (C2) (when T is the box operator).

If this point were correct, then it would be mysterious why anyone would have been tempted to blame the necessitation rule ($\Box$ Intro) for the paradox. Regardless, the point being made is not correct. After all, a derivation from $\sim\Box\mu$ to $\sim\mu$ is not sufficient for (C2) but rather for only a single instance of (C2). For (C2) to hold properly (when T is the box operator), it must be that for *any* wff $\phi$, there is a derivation from $\sim\Box\phi$ to $\sim\phi$. And it is plain that this general condition does not hold of the system T–; $\sim\Box\phi$ does not generally allow a derivation of $\sim\phi$.

It is correct, however, that $\sim\Box\mu$ allows one to derive $\sim\mu$. But that was made possible not because of a general condition like (C2) but rather because of the specific hypothesis occurring at line 1. The kind of contradiction enabled by conditions like (C1) and (C2) can instead be enabled by the hypothesis in conjunction with $\Box$ Elim, the first interchange rule, and a few truth-functional rules. (And without $\Box$ Intro.)

## *6. Closing*

In systems T and T–, Montague's paradox demonstrates that there is no sentence $\mu$ which is equivalent to its own negated necessitation. This in itself is not news. What has been underappreciated, however, is that rejecting the necessitation rule does not prevent us from demonstrating the result.[6]

## References

Boolos, G., Burgess, J., & Jeffrey, R. (2007). *Computability and Logic*, 5th edition. Cambridge: Cambridge University Press.

Dean, W. (2014). Montague's Paradox, Informal Provability, and Explicit Modal Logic. *Notre Dame Journal of Formal Logic* 55(2): 157–196.

Friedman, H. & Sheard, M. (1987). An Axiomatic Approach to Self-Referential Truth. *Annals of Pure and Applied Logic* 33: 1–21.

Garson, J. (2013). *Modal Logic for Philosophers*, 2nd edition. Cambridge: Cambridge University Press.

Halbach, V. (2009). On the Benefits of a Reduction of Modal Predicates to Modal Operators, in A. Hieke and H. Leitgeb, *Reduction—Abstraction—Analysis*, *proceedings of the 31st International Wittgenstein Symposium*. Kirchberg, Ontos Verlag, Heusenstamm bei Frankfurt a.M, pp. 323–333.

Kripke, S. (1959) 'A Completeness Theorem in Modal Logic,' *Journal of Symbolic Logic* 24: 1–14.

____. (1963). "Semantical Considerations on Modal Logic," *Acta Philosophica Fennica* 16: 83–94.

Montague, R. (1963). "Syntactical Treatments of Modality, with Corollaries on Reflexion Principles and Finite Axiomatizability," *Acta Philosophica Fennica* 16: 153–167.

Post, J. (1970). "The Possible Liar," *Noûs* 4(4): 405–409.

Quine, W.V.O. (1951). *Mathematical Logic*, revised edition. Cambridge, MA: Harvard University Press.

Robinson, R.M. (1950). "An Essentially Undecidable Axiom System," *Proceedings of the International Congress of Mathematics*: 729–730.

[6] Might an intuitionist have a different way to block the paradox? The truth-functional rules of T– are intuitionally acceptable, but perhaps the first interchange rule is not. For if '□' is interpreted as "it is constructively provable that…", the rule expresses the falsity that if $\phi$ is not constructively provable, then $\sim\phi \nvdash \bot$. But while this point is correct, it is unnecessary: Provability logics lack (T), and that is already enough to avoid the paradox.